\documentclass[11pt]{article}

\usepackage{multirow, multicol}   
\usepackage{booktabs} 
\usepackage{array} 
\usepackage{paralist} 
\usepackage{verbatim} 
\usepackage{subfig} 

\usepackage{amsmath, amssymb, enumerate}
\usepackage{pictexwd,dcpic}

\newtheorem{theorem}{Theorem}
\newtheorem{corollary}[theorem]{Corollary}
\newtheorem{remark}[theorem]{Remark}
\newtheorem{lemma}[theorem]{Lemma}
\newtheorem{definition}[theorem]{Definition}
\newtheorem{proposition}[theorem]{Proposition}
\newtheorem{example}[theorem]{Example}

\newcommand{\proof}{ {\sc Proof.\quad}}
\newcommand{\pend}{ \hfill $\square$ \\}

\numberwithin{equation}{section}  
\numberwithin{figure}{section}    
\numberwithin{table}{section}     
\numberwithin{theorem}{section}

\newcommand{\of}[1]{\ensuremath{\left( #1 \right)}}

\newcommand{\cb}[1]{\ensuremath{ \left\{ #1 \right\} }}
\newcommand{\sqb}[1]{\ensuremath{ \left[ #1 \right] }}

\newcommand{\st}{\,|\;}

\newcommand{\eps}{\ensuremath{\varepsilon}}

\newcommand{\vp}{\ensuremath{\varphi}}

\newcommand{\R}{\mathrm{I\negthinspace R}}
\newcommand{\OLR}{\overline{\mathrm{I\negthinspace R}}}
\newcommand{\N}{\mathrm{I\negthinspace N}}

\renewcommand{\P}{\ensuremath{\mathcal{P}}}

\newcommand{\dom}{{\rm dom \,}}

\newcommand{\Int}{{\rm int\,}}

\newcommand{\isum}{{+^{\negmedspace\centerdot\,}}}

\newcommand{\idif}{{-^{\negmedspace\centerdot\,}}}

\newcommand{\triup}{{\rm \vartriangle}}

\usepackage[
colorlinks=true, linkcolor=blue, citecolor=blue, urlcolor=blue, filecolor=blue
]{hyperref}

\begin{document}

\title{Variational inequalities characterizing weak minimality in set optimization}


\author{Giovanni P. Crespi\thanks{Department of Economics and Political Sciences, Loc. Grand Chemin 73-75,
		11020 Saint-Christophe, Aosta, Italy, \href{mailto: g.crespi@univda.it}{g.crespi@univda.it} }           
\and Matteo Rocca\thanks{Department of Economics, University of Insubria, via Ravasi, 2, 21100 Varese, Italy, 
		\href{mailto: matteo.rocca@uninsubria.it}{matteo.rocca@uninsubria.it} }
     \and Carola Schrage\thanks{Department of Economics and Political Sciences, Loc. Grand Chemin 73-75,
		11020 Saint-Christophe, Aosta, Italy, \href{mailto: carolaschrage@gmail.com}{carolaschrage@gmail.com}}        
}

\date{\small\today}

\maketitle

\begin{abstract}
We introduce the notion of weak minimizer in set optimization. Necessary and sufficient conditions in terms of scalarized variational inequalities of Stampacchia and Minty type, respectively, are proved.\\
As an application, we obtain necessary and sufficient optimality conditions for weak efficiency of vector optimization in infinite dimensional spaces. A Minty variational principle in this framework is proved as a corollary of our main result.
\end{abstract}


\section{Introduction}\label{sec:Intro}

Scalar variational inequalities (for short, VI) apply to study a wide range of practical problems, in particular
equilibrium and optimization problems,  see e.g. \cite{BaiocchiCapelo84}, \cite{KinderlehrerStempacchia80}.
Generalizations toward vector VI were initiated in  \cite{Giannessi80}; for recent results and survey on this field see e.g. \cite{AlHomidanAnsari2010}, \cite{Giannessi98}, \cite{GiannessiMastroeniPellegrini2000}, \cite{Komlosi99}, \cite{LiMastroeni2010}. Far less has been undertaken to extend those results to set-valued optimization, mainly because of a rather different approach to the classical optimization of set-valued maps.

In the scalar case, when the operator involved in a VI has a primitive function, we refer to the problem and a differentiable VI. This kind of VI is widely studied because of its relation to optimization problems. Under mild continuity assumptions, scalar Minty VI (MVI, \cite{LionsStampacchia67}, \cite{Minty67}) of differential type, provide a sufficient optimality condition to the primitive optimization problem (a result popularized as Minty variational principle), while scalar Stampacchia VI (SVI, \cite{Stampacchia64}) is only necessary. Assuming some convexity on the primitive function (or monotonicity of the derivative) both VIs are necessary and sufficient optimality conditions. In \cite{crespiJOTA}, under generalized differentiability assumptions, scalar Minty VI have been studied and it has been proved that the existence of a solution to such a problem implies some regularity property on the primitive optimization problem.

During the last decade, the Minty variational principle has been extended to the vector case. Since the seminal paper by Giannessi \cite{Giannessi98} the links between Minty variational inequalities and vector optimization problems were investigated. More recently, in \cite{CreGinRoc08}, \cite{YangYang04}, some generalizations of the vector principle have been proposed in conjuction with weak efficient solutions. In  \cite{Giannessi98}, \cite{YangYang04}, the case of a differentiable objective function $f$ with values in ${\mathbb R}^m$ and a Pareto ordering cone has been studied, proving a vector Minty variational principle for pseudoconvex functions. In \cite{CreGinRoc08} a similar result has been extended to the case of an arbitrary ordering cone and a nondifferentiable objective function. 

Optimization of set-valued functions has been a fast growing topic over the past decades. Since the first results by Corley \cite{Corley88}, \cite{corley87} and  Luc \cite{Luc89}, based on a vector optimization approach, several papers  provide optimality conditions. Nevertheless, the main approach to derivatives (and therefore to the core of a variational inequality) has been far distant form the basic differential quotient method adopted for scalar (and vector) problems. More recently, a new paradigm, known as set-optimization, has been proposed, compare \cite{HamelHabil05},\cite{HeydeLoehne11}, \cite{LoehneDiss}, \cite{Loehne11Book}. In this framework, the very concept of optimal solutions has been thought anew, together with operations among sets, now elements of a complete ordered conlinear space. This leads to overcome some drawbacks in previous attempt to provide variational inequalities for set-valued optimization problems (see e.g. \cite{CreGinRoc2010}).

In this paper, we present the notion of weak minimality for set-optimization, motivated by its relation with standard weak efficiency in vector optimization. Using scalarization techniques, we define Minty and Stampacchia type differential variational inequalities corresponding to the primitive set--optimization problem. We prove, under Hausdorff continuity and suitable pseudoconvexity assumptions that the solutions of the Minty type inequality are weak minimizers of the primitive set-optimization problem. Under slightly weaker assumptions, a weak minimizer of the set-optimization problem solves the Stampacchia differential variational inequality. Under convexity assumptions on the scalarizations, the reverse implications has been proven in \cite{CrespiSchrage13b}. As an application, we obtain new results for vector optimization in infinite dimensional spaces that are comparable to those in \cite{CreGinRoc08} for the finite dimensional case.
A similar approach, but for minimizers rather then weak minimizers can be found in \cite{CrespiSchrage13a} and in finite dimensions in \cite{pilecka14_Preprint}.

The paper is organized as follows. Section 2 is devoted to some preliminary results, notation and definitions that are used throughout the paper. We introduce the general setting for the problem and the scalarization technique that is used to prove the main results. The concept of minimality is also introduced and commented together with the main properties that constitutes assumptions for the results provided in the next sections. The main results are proved in Section 3, both for Minty and Stampacchia variational inequality. Each theorem of Section 3 is matched with a corollary that proves vector optimization result as a special case. Finally, in Section 4 we draw some conclusions and provide some insight on the complete ring that we are developing in conjunction with previous results. Indeed, the paper fits in a line of research that includes other results and ideas. Some of them are part of other papers, that are currently submitted or not yet published while this paper is being drafted. For the readers convenience we add an appendix with the proofs of those results that are used in this paper, originally proved in other papers yet to appear.

\section{Preliminaries}\label{sec:Preliminaries}

Let $X,Y$ be locally convex Hausdorff spaces with topological duals $X^*, Y^*$. The set $\mathcal U_Y(0)$ ($\mathcal U_X(0)$) is a $0$-neighborhood base of $Y$ (of $X$) consisting of balanced convex open sets.
 $Y$ is pre-ordered by a closed convex cone $C$ with nonempty interior $\Int C\neq\emptyset$ and $C\neq Y$. We denote by $\P(Y)$ the set of all subsets of $Y$. 
For all $A,B\in\P(Y)$ we set
\[
A< B\quad\Leftrightarrow\quad B\subseteq A+\Int C,
\]
compare i.e \cite[Definition 3.2]{maeda2010}.

\begin{lemma}
For all $A\in\P(Y)$ it holds $ A+\Int C=\Int (A+C)$.
\end{lemma}
\proof
As the sum of an open set and an arbitrary set is always open, it is only left to prove that $a\in \Int(A+C)$ implies $a\in A+\Int C$. Thus, assume $a+U\subseteq A+C$ is satisfied for some $U\in\mathcal U_Y(0)$.
Let $e\in \Int C$ and $n\in \N$ be such that $-\frac{1}{n}e\in U$, then
\[
a\in A+\frac{1}{n}e+C\subseteq A+(\Int C+C)\subseteq A+\Int C.
\]
\pend

Since, $C+\Int C=\Int C$, for all $A, B\in\P(Y)$ it holds
\[
A<B\quad\Leftrightarrow \of{A+C} < \of{B+C}.
\]
Strict inequality between $A,B\in\P(Y)$ is denoted by
\[
A\ll B\quad\Leftrightarrow\quad  \exists U\in\mathcal U_Y(0):\;B+U\subseteq A+C.
\]

It is an easy task to prove that for all $A, B\in\P(Y)$ it holds
\[
A\ll B\quad\Leftrightarrow\quad(A+C)\ll(B+C).
\]

\begin{proposition}\label{prop:minim_connect}
Let $A, B\in\P(Y)$, then
$A\ll B$ implies $A<B$. If additionally $B$ is compact, then the reverse implication holds as well.
\end{proposition}
\proof
The first implication is immediate. Assume by contradiction $B$ is compact, $B\subseteq \Int\of{A+ C}$
 and for all $U\in\mathcal U_Y(0)$ there exists a $b_U\in B$, $u\in U$ with  $(b_U+u)\notin \of{A+C}$.
Let $I$ be a nonempty index set, $\cb{U_i}_{i\in I}\subseteq \mathcal U_Y(0)$ be given with $U_i\subseteq U_j$ and $U_i\neq U_j$ whenever $j< i$ and
let $\cb{b_i}_{i\in I}\subseteq B$ be given with $b_i =b_{U_i}$.
Compactness of $B$ implies the existence of a convergent subnet, hence without loss  assume $b_i\to b_0\in B$.
By assumption, there exists  $U_0\in\mathcal U_Y(0)$ with $b_0+U_0\subseteq A+C$ and for $i\in I$ large enough, $b_i+U_i\subseteq b_0+U_0$. But this is a contradiction, as we assumed $b_i+U_i\nsubseteq A+C$ and $b_0+U_0\subseteq A+C$.
\pend

The reverse implication is not true in general.
\begin{example}
Let $B=\cb{(x,\frac{1}{x})\st x>0}$ and $A=C=\R^2_+$ in $\R^2$. Clearly, $B\subseteq A+\Int C$ while $(B+U)\setminus A\neq\emptyset$ is true for all $U\in\mathcal U_{\R^2}(0)$.
\end{example}

The positive dual cone of $C$ is the set $C^+=\cb{y^*\in Y^*\st \inf\limits_{c\in C}\cb{y^*(c)}\geq 0}$.\\
Since $\Int C\neq\emptyset$ is assumed, $C^+$ possesses a weak$^*$-compact base $W^*$, see e.g. \cite[Theorem 1.5.1]{Aubin71}.

\begin{remark}\label{rem:boundedValuesonU}\cite[Remark 3.32]{HeydeSchrage11R}
Under our assumptions, the following inequalities are satisfied
\begin{align*}
&\forall y\in Y:\inf\limits_{w^*\in W^*}w^*(y)>-\infty;\\
&\forall U\in\mathcal U_Y(0):\quad \sup\limits_{w^*\in W^*}\inf\limits_{u\in U}w^*(u)<0.
\end{align*}
\end{remark}

\begin{remark}\cite[Proposition 2.11]{CrespiSchrage13b}
  Let $A,B\in\P(Y)$ be given such that $B+C$ is convex. Then
  $A\ll B$ implies
  \begin{align*}
  \forall w^*\in W^*:\quad \inf\limits_{a\in A}w^*(a)=-\infty\;\vee\;\inf\limits_{a\in A}w^*(a)<\inf\limits_{b\in B}w^*(b),
  \end{align*}
  which in turn implies $A<B$ if also $A+C$ is convex. Under the additional assumption of compactness of $B$, equivalence is proven by Proposition \ref{prop:minim_connect}.
  Moreover, under compactness, $\inf\limits_{b\in B}w^*(b)$ is attained at some $b\in B$ for all $w^*\in W^*$ whenever $B$ is nonempty, in which case the value especially is finite.
\end{remark}

For any function $F:X\to\P(Y)$ and $y^*\in C^+\setminus\cb{0}$, we define the scalarization of $F$ w.r.t. $y^*$ as
\[
\vp^\triup_{F,y^*}:X\to\OLR,\; \vp^\triup_{F,y^*}(x)=\inf\cb{y^*(y)\st y\in F(x)}.
\]

Since
\[
\vp^\triup_{F,y^*}(x)=+\infty\quad\Leftrightarrow\quad   F(x)= \emptyset,
\]
the effective domain of any scalarization
\[
\dom\vp^\triup_{F,y^*}=\cb{x\in X\st \vp^\triup_{F,y^*}(x)\neq+\infty}
\]
 coincides with that of the set valued function,
\[
\dom F=\cb{x\in X\st F(x)\neq \emptyset}.
\]

It easily follows that when $F(x)$ is compact for all $x\in X$, then $\dom F\neq\emptyset$, if and only if $\vp^\triup_{F,w^*}:X\to\OLR$ is proper for all $w^*\in W^*$, i.e $\dom \vp^\triup_{F,w^*}\neq\emptyset$ and $\vp^\triup_{F,w^*}$ does not attain the value $-\infty$.

\begin{definition}\label{def:weak_minimizer}\cite{CrespiSchrage13b}\cite{HernandezMarin05}
Let $F:X\to\P(Y)$ be given.
 Then $x_0\in\dom F$ is called {\em weak-l}, {\em scalarized weak} or {\em weak  minimizer} of $F$, if either $F(x_0)+C=Y$, or
\begin{align}\label{eq:w-l-Min}
\tag{$w$-$l$-$Min$}
\forall x\in X:\quad F(x)\nless F(x_0);
\end{align}
\begin{align}\label{eq:w-sc-Min}
\tag{$w$-$sc$-$Min$}
\forall x\in X:\exists w^*\in W^*: \quad \vp^\triup_{F,w^*}(x_0)\leq \vp^\triup_{F,w^*}(x)\neq-\infty;
\end{align}
\begin{align}\label{eq:w-Min}
\tag{$w$-$Min$}
\forall x\in X:\quad F(x)\not\ll F(x_0).
\end{align}
\end{definition}

Without further assumptions, \eqref{eq:w-l-Min} implies \eqref{eq:w-Min}. The equivalence holds if $F(x_0)$ is compact. If $F(x_0)+C$ is convex, then
\[
\eqref{eq:w-l-Min} \Rightarrow \eqref{eq:w-sc-Min} \Rightarrow \eqref{eq:w-Min}
\]
Moreover, if  $F(x_0)$ is compact, then the three notions are equivalent. 

\begin{remark}
If $F:X\to Y$ is a vector valued function, then Definition \ref{def:weak_minimizer} reduces to the notion of weak efficiency for  vector optimization.
\end{remark}

\begin{definition}
A set valued function $F:X\to\P(Y)$ is called $C$-convex, when for all $x_1,x_2\in X$ and all $t\in\sqb{0,1}$ it holds
\[
t F(x_1)+(1-t)F(y_2)\subseteq F(tx_1+(1-t)x_2)+C.
\]
\end{definition}

\begin{remark}\cite[Corollary 1.69]{Loehne11Book}
A set valued function $F:X\to\P(Y)$ is $C$-convex, 
if and only if the functions
$\vp^\triup_{F,y^*}:X\to\OLR$
are convex for all $y^*\in C^+\setminus\cb{0}$.
In this case,  $F(x)+C$ is a convex subset of $Y$ for all $x\in X$.
\end{remark}

\begin{definition}\cite[Definition 4.4]{CrespiHamelSchrage13W}
A scalar function $\vp:X\to\OLR$ is semistrictly quasiconvex when
\[
\forall a,b\in\dom\vp,\,\forall t\in\of{0,1}:\quad \vp(a) \neq \vp(b)\;\Rightarrow\; \vp(a+t(b-a)) < \max\cb{\vp(a),\vp(b)}.
\]
It is (lower Dini) pseudoconvex, when $\vp(a) < \vp(b)$ implies
\[
\vp^\downarrow(b,a-b)=\liminf\limits_{t\downarrow 0}\frac{1}{t}\inf\cb{\rho\in\R\st \vp(b+t(a-b))\leq \vp(b)+\rho} < 0
\]
and (lower Dini) pseudoconcave, when $\vp(a)>\vp(b)$ implies
\[
\vp^\downarrow(b,a-b) > 0.
\]
\end{definition}

The following property is used in the sequel.

\begin{proposition}\cite[Proposition 4.13]{CrespiHamelSchrage13W}\label{prop:char_semstr_qconvex}
Let $\vp :\R \to \OLR$ be l.s.c. and semistrictly quasiconvex with $\dom \vp\subseteq \sqb{0,1}$. Then there exist $s_0 \leq t_0 \in \sqb{0,1}$ such that $\vp$ is strictly decreasing on $\of{0, s_0}$, strictly increasing on $\of{t_0,1}$ and constantly equal to $\inf\cb{\vp(x)\st x\in X}$ on $\sqb{s_0, t_0}$.
\end{proposition}

Eventually, to prove our main results we can weaken some assumption to hold only on restriction along rays. \\
Given a single valued function $\vp:X\to\OLR$ and two points $x_0,x\in X$, we introduce the restriction of $\vp$ along the ray with extreme points $x$ and $x_0$ as $\vp_{x_0,x}:\R\to\OLR$ defined by
\[
\vp_{x_0,x}(t)=\begin{cases}
\vp(x_0+t(x-x_0)),&\text{ if } 0\leq t\leq 1;\\
+\infty.&\text{ elsewhere.}
\end{cases}
\]

Then $\vp:X\to\OLR$ is called radially semistrictly quasiconvex,  pseudoconvex or pseudoconcave at $x_0$, if $\vp_{x_0,x}$ is semistrictly quasiconvex, pseudoconvex or pseudoconcave for all $x\in X$.
Likewise, $\vp:X\to\OLR$ is called radially lower semicontinuous (l.s.c.) at $x_0$, if $\vp_{x_0,x}$ is l.s.c. for all $x\in X$.

In \cite[Proposition 4.14]{CrespiHamelSchrage13W} it is proven that if
 $\dom \vp$ is star shaped at $x_0$, i.e. $x_0,x\in\dom \vp$ implies the whole interval $\cb{x_0+t(x-x_0)\st 0\leq t\leq 1}$ is a subset of $\dom \vp$,  and  $\vp$ is radially pseudoconvex and l.s.c. at $x_0$, then it is radially semistrictly quasiconvex at $x_0$.
 A scalar convex, strictly monotone function on the real line  is  pseudoconvex and pseudoconcave.
Therefore, if $F:X\to\P(Y)$ is $C$-convex, $y^*\in Y^*$ and $\of{\vp^\triup_{F,y^*}}_{x_0,x}:\R\to\OLR$ is strictly decreasing on the interval $\sqb{0,s_0}\subseteq \dom F$, then
\[
\forall s\in\of{0,s_0}:\quad \of{\vp^\triup_{F,y^*}}^\downarrow_{x_0,x}(s,-1)>0.
\]

Finally, when dealing with Minty-type variational principle, some continuity is needed. For set--valued functions, we consider the following notions.
\begin{definition} (compare \cite[Proposition 1.5.2]{AubinCellina84})
A set valued function $F:X\to\P(Y)$ is  upper Hausdorff continuous in $x_0$, if
for all $U\in\mathcal U_Y(0)$ there exists a $0$-neighborhood $W\subseteq X$ with
\begin{align*}
\forall x\in x_0+V:\quad  F(x)\subseteq F(x_0)+U.
\end{align*}
\end{definition}

\begin{definition}\cite[Definition 1.5.2]{Luc89}
A set $\Psi=\cb{\vp_i:X\to\OLR\st i\in I}$ is lower equicontinuous in $x_0\in \bigcap\limits_{i\in I}\dom\vp_i$, if
\[
\forall \eps>0,\,\exists V\in\mathcal U_X(0),\,\forall x\in x_0+V,\,\forall i\in I:\quad
\vp_i(x_0)\leq \vp_i(x)+\eps
\]
\end{definition}

Hausdorff continuity of $F$ is related to lower equicontinuity  of its scalarizations.

\begin{proposition}\label{prop:Hausdorff_equ_Equicont}
If $F:X\to\P(Y)$ is upper Hausdorff continuous in $x_0\in \dom F$, then
$\Psi=\cb{\vp^\triup_{F,w^*}:X\to\OLR\st w^*\in W^*}$ is lower equicontinuous in $x_0\in X$.

If $F(x_0)+C$ is nonempty, convex and $\Psi$ is lower equicontinuous in $x_0$, then $F^C:X\to \P(Y)$ is upper Hausdorff continuous in $x_0$, defining 
\begin{align}\label{eq:sv_ext}
\forall x\in X:\quad
F^C(x)=
\begin{cases}
F(x)+C,&\text{ if } F(x)\neq\emptyset;\\
\emptyset,&\text{ elsewhere}.
\end{cases}
\end{align}
\end{proposition}
\proof
First assume $F$ and thus $F^C$ is upper Hausdorff continuous in $x_0\in\dom F$, $-e\in\Int C$. By Remark \ref{rem:boundedValuesonU}, without loss of generality, we can assume
\[
\inf\limits_{w^*\in W^*}w^*(e)=-1
\]
and $(\eps e+C)$ contains a neighborhood $U\in \mathcal  U_Y(0)$ for all $\eps>0$, implying the existence of $V\in\mathcal U_X(0)$ such that for all $x\in x_0+V$ it holds
\[
\forall w^*\in W^*:\quad \vp^\triup_{F,w^*}(x)\geq \vp^\triup_{F,w^*}(x_0)-\eps,
\]
which is lower equicontinuity of $\Psi$ in $x_0$.\\
On the other hand, assume $F(x_0)+C$ is a convex set and
\[
\exists U\in\mathcal U_Y(0),\, \forall V\in \mathcal U_X(0),\, \exists x\in x_0+V:\quad F(x)\nsubseteq F(x_0)+C+U.
\]
Especially, for all $V\in \mathcal U_X(0)$ there exists $x\in x_0+V$ such that, by a separation argument,
\[
\exists w^*\in W^*:\quad \vp^\triup_{F,w^*}(x)\leq \vp^\triup_{F,w^*}(x_0)+\inf\cb{w^*(u)\st u\in U}\in\R.
\]
But, by Remark \ref{rem:boundedValuesonU},
\[
\sup\limits_{w^*\in W^*}\inf\cb{w^*(u)\st u\in U}=-\mu<0,
\]
hence there exists $\mu>0$ such that
\[
\forall V\in\mathcal U_X(0),\,\exists x\in x_0+V,\,\exists w^*\in W^*:\quad \vp^\triup_{F,w^*}(x)< \vp^\triup_{F,w^*}(x_0)-\frac{1}{2}\mu,
\]
contradicting lower equicontinuity of $\Psi$ in $x_0$.
\pend

\begin{remark}
In the literature, there are many competing concepts of continuity notions for set valued functions. As an example,
$C$-upper continuity of $F$ as defined in \cite[Definition 2.2]{KuwanoTanaka12} implies 
upper Hausdorff continuity of $F^C$, while  upper Hausdorff continuity of $F$ implies local Lipschitz continuity of $F$, as defined in \cite[Definition 2.3]{KuwanoTanaka12}.
For a closer study of different continuity concepts and their correlation, compare \cite{HeydeSchrage11R}.
\end{remark}

We say that $F$ is radially upper Hausdorff continous w.r.t $x_0$, if  $F_{x_0,x}$ is upper Hausdorff continous in $\sqb{0,1}$ for all  $x\in X$,
setting $F_{x_0,x}:\R\to\P(Y)$ as
\[
F_{x_0,x}(t)=
\begin{cases}
F(x_0+t(x-x_0)),&\text{ if } 0\leq t\leq 1;\\
\emptyset,& \text{ elsewhere.}
\end{cases}
\]
Moreover
$$
\of{\vp^\triup_{F,w^*}}_{x_0,x}(t)=\vp^\triup_{F_{x_0,x},w^*}(t)
$$  
is true for all $w^*\in W^*$, $x\in X$  and all $t\in\R$.
Thus if $F$ is radially upper Hausdorff continous w.r.t $x_0$, then for all $x\in X$ the set
\[
\Psi_{x_0,x}=\cb{\of{\vp^\triup_{F,w^*}}_{x_0,x}:\R\to\OLR\st w^*\in W^*}
\]
is lower equicontinuous in $t$ for all $0\leq t\leq 1$. If additionally $F(x)=F(x)+C$ for all $x\in X$ and the images are convex, then the equivalence holds.

\begin{lemma}\label{lem:continuity}
Let $F:X\to\P(Y)$ be such that there exists $x_0\in X$ with $F(x_0)$ compact and
let $\cb{w^*_i}_{i\in I}\subseteq W^*$, $I$ a nonempty index set and $w^*_i\to w^*_0$ in the weak$^*$ topology.
Then
\[
\lim\limits_{w^*_i\to w^*_0}\vp^\triup_{F,w_i^*}(x_0)= \vp^\triup_{F,w^*_0}(x_0)
\]
\end{lemma}
\proof
Compactness of $F(x_0)$ implies that each $w_i^*\in W^*$ has a supporting point $z_i\in F(x_0)$ to $F(x_0)+C$. Without loss of generality, assume $z_i\to z_0\in F(x_0)$ is satisfied.
Weak$^*$ compactness of $W^*$ implies that
\[
\forall \mu>0,\,\exists i_\mu,\,\forall i>i_\mu: w^*_i(z_0)-\mu\leq w^*_i(z_i)=\vp^\triup_{F,w^*_i}(x_0)\in\R.
\]
As the support function of a set is weak$^*$ l.s.c. and convex, $w^*\mapsto\vp^\triup_{F,w^*}(x_0)$ is concave and weak$^*$ upper semicontinuous and it holds
\[
\limsup\limits_{w^*_i\to w^*_0}\vp^\triup_{F,w^*_i}(x_0)\leq \vp^\triup_{F,w^*_0}(x_0)
\leq w^*_0(z_0)\leq \liminf\limits_{w^*_i\to w^*_0}\vp^\triup_{F,w^*_i}(x_0),
\]
Thus $\lim\limits_{w^*_i\to w^*_0}\vp^\triup_{F,w^*_i}(x_0)=\vp^\triup_{F,w^*_0}(x_0)$ and $z_0$ is a supporting point of $F(x_0)+C$ to $w^*_0$.
\pend

\begin{lemma}\label{lem:combined_lsc}
Let $F:X\to\P(Y)$ be such that $\Psi=\cb{\vp^\triup_{F,w^*}:X\to\OLR\st w^*\in W^*}$ is lower equicontinuous in $x_0\in \dom F$ and
$F(x_0)$ is compact.
Let $\cb{w^*_i}_{i\in I}\subseteq W^*$, $I$ a nonempty index set and $w^*_i\to w^*_0$ in the weak$^*$ topology, then
\[
\liminf\limits_{\substack{x\to x_0\\w^*_i\to w^*_0}}\vp^\triup_{F,w^*}(x)\geq \vp^\triup_{F,w^*_0}(x_0).
\]
\end{lemma}
\proof
As by assumption $F(x_0)$ is compact, $\vp^\triup_{F,w^*}(x_0)\in\R$ is true for all $w^*\in W^*$ and for all $n\in\N$ there exists $V\in\mathcal U_X(0)$ such that
\[
\forall w^*\in W^*,\,\forall x\in x_0+V:\quad \vp^\triup_{F,w^*}(x_0)\leq \vp^\triup_{F,w^*}(x)+\frac{1}{n}.
\]
But by Lemma \ref{lem:continuity}, eventually $\vp^\triup_{F,w^*_0}(x_0)\leq \vp^\triup_{F,w^*}(x_0)+\delta$ for all $\delta>0$ as $w^*$ converges to $w^*_0$. Thus
\[
\liminf\limits_{\substack{x\to x_0\\w^*\to w^*_0}}\vp^\triup_{F,w^*}(x)\geq \vp^\triup_{F,w^*_0}(x_0)
\]
is true.
\pend


\section{Main results}\label{sec:Main}

\subsection{Minty variational principle}\label{subsec:Minty}\

Minty variational principle (see e.g. \cite{Giannessi98}, \cite{Minty67}) provides a sufficient optimality condition in terms of a variational inequality under mild continuity assumptions.\\
Recent results (see e.g. \cite{crespiJOTA}, \cite{CreGinRoc08}, \cite{YangYang04}) have formalized the variational inequality by means of a generalized dini--type directional derivative.\\
Using the scalarizations $\vp^\triup_{F,W^*}$ we prove sufficient optimality condition for weak--minimizers of a set--valued function $F$ under Hausdorff continuity assumption. Since Theorem \ref{thm:MainResult} is stated through a scalarized Minty variational inequality (\ref{eq:mvi_scalar}), we can interpret it as a Minty variational principle for vector optimization.

\begin{theorem}\label{thm:MainResult}
Let $F:X\to\P(Y)$ be radially upper Hausdorff continuous at $x_0\in\dom F$, $\dom F$ be star shaped at $x_0$ and $F(x_0)$ be compact.
Moreover, assume $\vp^\triup_{F,w^*}$ is proper  and radially pseudoconvex and pseudoconcave w.r.t. $x_0$
for all $w^*\in W^*$.
If
\begin{align}\label{eq:mvi_scalar}\tag{$mvi$}
\forall x\in X,\, \exists w^*\in W^*:\quad \of{\vp^\triup_{F,w^*}}^\downarrow(x,x_0-x)\leq 0
\end{align}
is satisfied, then $x_0$ is a weak minimizer of $F$.
\end{theorem}
\proof
Radially upper Hausdorff continuity at $x_0$ implies $\vp^\triup_{F,w^*}$ is l.s.c. on the interval $\cb{x_0+t(x-x_0)\st 0\leq t\leq 1}$ for all $w^*\in W^*$ and all $x\in X$.
Next, assume to the contrary that $F(x_0)\neq Y$ and there exists a $x\in X$ such that
\[
F(x_0)+U\subseteq F(x)+C
\]
is true for some neighborhood $U\in\mathcal U_Y(0)$.
By Remark \ref{rem:boundedValuesonU}
\[
\sup\limits_{w^*\in W^*}\inf\limits_{u\in U}w^*(u)=-\mu<0,
\]
thus,
properness of the scalarizations $\vp^\triup_{F,w^*}$ implies
\[
\forall w^*\in W^*:\quad
-\infty\neq \vp^\triup_{F,w^*}(x)-\vp^\triup_{F,w^*}(x_0)\leq-\mu<0
\]

Pseudoconvexity and lower semicontinuity imply semistrict quasiconvexity, see e.g. \cite[Proposition 4.13]{CrespiHamelSchrage13W}. Thus there exists $0<t_{w^*}\leq 1$ such that $\vp^\triup_{F,w^*}$ is strictly decreasing on the interval $\cb{x_0+t(x-x_0)\st t\in\sqb{0,t_{w^*}}}$ as $t$ converges to $t_{w^*}$ and $\vp^\triup_{F,w^*}$ is constant on the interval  $\cb{x_0+t(x-x_0)\st t\in\sqb{t_{w^*},1}}$.

For all $w^*\in W^*$, we define the function $\Phi_{F,w^*}:\R\to\OLR$
\[
\Phi_{F,w^*}(t)=
\begin{cases}
\vp^\triup_{F,w^*}(x_0+t(x-x_0))- \vp^\triup_{F,w^*}(x_0),&\text{ if } 0\leq t\leq 1;\\
+\infty,&\text{ elsewhere.}
\end{cases}
\]
that is pseudoconvex, l.s.c. in $\sqb{0,1}$, attains a global minimum in $t_{w^*}$ and
\[
\Phi_{F,w^*}(t_{w^*})\leq \sup\limits_{w^*\in W^*} \Phi_{F,w^*}(1)\leq -\mu<0.
\]

Next, assume
\[
t_0=\inf\cb{t_{w^*}\st w^*\in W^*}=0.
\]
Especially we can find a convergent net $\cb{w^*_i}_{i\in I}\in W^*$, $w^*_i\to w^*_0\in W^*$.
As
$F(x_0)$ is compact, by Lemma \ref{lem:combined_lsc} and Proposition \ref{prop:Hausdorff_equ_Equicont} upper Hausdorff continuity in $x_0$ implies
\[
-\infty\neq \vp^\triup_{F,w_0^*}(x_0)\leq \liminf\limits_{i\in I}\vp^\triup_{F,y_i^*}(x_0+t_{w^*_i}(x-x_0))
\]
contradicting
\[
\liminf\limits_{i\in I}\Phi_{F,y_i^*}(t_{w^*_i})\leq -\mu<0.
\]
Hence $t_0>0$ and for all $w^*\in W^*$, $\vp^\triup_{F,w^*}$ is strictly decreasing on the interval  $\cb{x_0,x_0+t(x-x_0)\st t\in\sqb{0,t_0}}$.
Setting $\bar x= x_0+\frac{1}{2}t_0(x-x_0)$, then
pseudoconcavity implies
\[
\forall w^*\in W^*:\quad \of{\vp^\triup_{F,w^*}}^\downarrow(\bar x,(x_0-\bar x)> 0,
\]
contradicting \eqref{eq:mvi_scalar}, proving the statement.
\pend

\begin{remark}
The assumption $F\left(x_0\right)$ compact implies that weak minimizers coincides with weak-l and scalarized weak ones. Therefore Theorem \ref{thm:MainResult} provides a sufficient condition for any notion of minimality in Definition \ref{def:weak_minimizer}
\end{remark}

In Theorem \ref{thm:MainResult}, the assumption  $\vp^\triup_{F,w^*}$ radially pseudoconvex and pseudoconcave w.r.t. $x_0$ for all $w^*\in W^*$ can be replaced by (radial) $C$-convexity of $F$.

\begin{corollary}
Let  $F:X\to\P(Y)$ be compact valued, $x_0\in\dom F$ and let $F$ be radially upper  Hausdorff continuous w.r.t. $x_0$. If  $F$ is $C$-convex, then \eqref{eq:mvi_scalar} implies $x_0$ is a weak (weak-l, scalarized weak) minimizer of $F$.
\end{corollary}
\proof
Compactness of $F(x)$ implies $\vp^\triup_{F,w^*}(x)\in\R$ for all $w^*\in W^*$, if $x\in\dom F$.
$C$-convexity and radial upper Hausdorff continuity of $F$ (w.r.t. $x_0$) imply radial pseudoconvexity, radial semistrict quasiconvexity of $\vp^\triup_{F,w^*}$  w.r.t. $x_0$ and radial lower equisemicontinuity of the scalarizations $\vp_{F,w^*}$.
Moreover,
\[
\forall s\in\of{0,s_0}:\quad \of{\vp^\triup_{F,w^*}}^\downarrow_{x_0,x}(s,-1)>0
\]
whenever $\of{\vp^\triup_{F,w^*}}_{x_0,x}$ is strictly decreasing on the interval $\sqb{0,s_0}$ is true, if $F$ is $C$-convex.
The domain of $F$ is convex, thus star shaped at $x_0\in \dom F$.
Applying the same arguments as in the proof of Theorem \ref{thm:MainResult}, \eqref{eq:mvi_scalar} implies $x_0$ is a weak minimizer of $F$.
But compactness of the images of $F$ combined with $C$-convexity of $F$ implies that in this case, $x_0$ is a weak-l and a scalarized weak minimizer of $F$ as well.
\pend



As an application, we can prove as a special case a result on vector optimization.

\begin{corollary}\label{3.4}
Let $S$ be star shaped at $x_0\in S\subseteq X$,  $F:S\to Y$ radially $C$-continuous w.r.t $x_0$, i.e. for all $x\in S$ it holds
\[
\forall U\in\mathcal U_Y(0),\,\exists t_U\in\of{0,1},\,\forall t\in \of{0,t_U}: F(x_0+t(x-x_0))\in F(x_0)+C+U.
\]
If $w^*\circ F:X\to\OLR$ is either radially pseudoconvex and radially pseudoconcave or radially convex w.r.t. $x_0$ for all $w^*\in W^*$,
 then \eqref{eq:mvi_scalar} implies $(F(x_0)-\Int C)\cap\cb{F(x)\st x\in S}=\emptyset$, i.e. $x_0$ is a weakly efficient minimizer  of   $F$.
\end{corollary}

The previous result can be compared with \cite[Theorem 3.7]{CreGinRoc08}. However, the older result is stated for functions with pseudoconvex scalarizations into image spaces with finite dimension. So, while we loose generality as we need a stronger convexity assumption, we allow for the more general setting on infinite dimensional spaces. Therefore Corollary \ref{3.4} is a new result also for vector optimization. 


\subsection{Stampacchia variational principle}\label{subsec:Stamp}

Necessary conditions of variational type can be proved through a slightly different type of variational inequality. Namely, the directional derivative would be evaluated at $x_0$ instead of $x$.

\begin{theorem}\label{3.5}
Let $F:X\to\P(Y)$  with $\emptyset\neq F(x_0)$ compact be such that all scalarizations $\vp^\triup_{F,w^*}:X\to\OLR$ with $w^*\in W^*$ be radially l.s.c. and  semistrictly quasiconvex  w.r.t. $x_0$. If $x_0$ is a weak minimizer, then we have
\begin{align}\label{eq:svi_scalar}\tag{$svi$}
\forall x\in \dom F,\, \exists w^*\in W^*:\quad \of{\vp^\triup_{F,w^*}}^\downarrow(x_0,x-x_0)\geq 0.
\end{align}
\end{theorem}
\proof
Assume to the contrary that
\begin{align*}
\exists x\in \dom F,\, \forall w^*\in W^*:\quad \of{\vp^\triup_{F,w^*}}^\downarrow(x_0,x-x_0)< 0.
\end{align*}
This implies that
\[
t_{w^*}=\sup\cb{t\in\sqb{0,1}\st \vp^\triup_{F,w^*}(x_0+t(x-x_0))<\vp^\triup_{F,w^*}(x_0)}>0
\]
for all $w^*\in W^*$. We set $t_0=\inf\cb{t_{w^*}\st w^*\in W^*}$ and without loss of generality assume $w^*\to w^*_0\in W^*$ as $t_{w^*}\to t_0$.

If $t_0>0$ is true, then by Proposition \ref{prop:char_semstr_qconvex}
\begin{align*}
\forall w^*\in W^*:\quad \vp^\triup_{F,w^*}(x_0+\frac{1}{2}t_0(x-x_0))<\vp^\triup_{F,w^*}(x_0),
\end{align*}
thus $x_0$ is does not satisfy \eqref{eq:w-sc-Min}.

On the other hand, assume $t_0=0$. As $\of{\vp^\triup_{F,w^*_0}}^\downarrow(x_0,x-x_0)<0$ is assumed, $t_{w^*_0}>0$ is true and applying Lemma \ref{lem:continuity} an Proposition \ref{prop:char_semstr_qconvex} we conclude
\[
\limsup\limits_{w^*\to w^*_0}\vp^\triup_{F,w^*}(x_0+\frac{1}{2}t_{w^*_0}(x-x_0))
\leq\vp^\triup_{F,w^*_0}(x_0+\frac{1}{2}t_{w^*_0}(x-x_0))
<\vp^\triup_{F,w^*_0}(x_0)
\]
But as $t_{w^*}\to 0$ is assumed,  $\frac{1}{2}t_{w^*_0}>t_{w^*}$ is true eventually, thus by semistrict quasiconvexity
\[
\vp^\triup_{F,w^*}(x_0)\leq\vp^\triup_{F,w^*}(x_0+\frac{1}{2}t_{w^*_0}(x-x_0))
\]
eventually and  by Lemma \ref{lem:continuity}
\[
\vp^\triup_{F,w^*_0}(x_0) =
\limsup\limits_{w^*\to w^*_0}\vp^\triup_{F,w^*}(x_0)\leq
\limsup\limits_{w^*\to w^*_0}\vp^\triup_{F,w^*}(x_0+\frac{1}{2}t_{w^*_0}(x-x_0))<\vp^\triup_{F,w^*_0}(x_0),
\]
a contradiction.
\pend

\begin{remark}
Since $F\left(x_0\right)$ is compact, weak minimizers coincide with weak-l and scalarized weak ones.
\end{remark}

We can apply Theorem \ref{3.5} to vector optimization in order to prove a necessary condition for weak efficiency.

\begin{corollary}
Let $F:S\subseteq X\to Y$ be radially $C$-convex and radially $C$-continuous w.r.t. $x_0\in\dom F$. If $x_0$ is a weakly efficient solution to $F$, then
\begin{align*}
\forall x\in \dom F,\, \exists w^*\in W^*:\quad \of{w^*\circ F}^\downarrow(x_0,x-x_0)\geq 0
\end{align*}
is satisfied.
\end{corollary}

\section{Conclusion}\label{sec:Concl}

In order to develop a complete loop between set optimization and variational inequalities, we need to pay some attention to infinite values of the scalarization. The assumption of compact images here excludes this possibility, that has been considered in \cite{CrespiSchrage13b}.
\\
For the sake of completeness we quote the result needed to complete the picture in Corollary \ref{cor:conclusion}.
\begin{proposition}\cite[Proposition 4.5]{CrespiSchrage13b}\label{prop:Stamp_Min2}
Let $F:X\to \P(Y)$ be a $C$-convex function, $x_0\in \dom F$.
If $x_0$ solves the scalarized Stampacchia variational inequality
\begin{align}\label{eq:svi2}
F(x_0)+C=Y\quad \vee\quad\forall x\in X\, \exists w^*\in W^*:\quad \of{\vp^\triup_{F,w^*}}^\downarrow(x_0,x-x_0)\geq 0,
\end{align}
then it is a scalarized weak minimizer.
\end{proposition}
The elements $x\in X\setminus\dom F$ can easily be discussed, as  $x_t=x_0+t(x-x_0)\in\dom F$ implies
\[
\of{\vp^\triup_{F,w^*}}^\downarrow(x_0,x-x_0)=\of{\vp^\triup_{F,w^*}}^\downarrow(x_0,x_t-x_0),
\]
for all $t>0$, while if $x_t\notin\dom F$ for all $t>0$ implies
\[
\of{\vp^\triup_{F,w^*}}^\downarrow(x_0,x-x_0)=+\infty,
\]
as by convention in \cite{CrespiSchrage13b} $(+\infty)\idif r=+\infty$ is true for all $r<+\infty$.

\begin{proposition}\cite[Lemma 4.9]{CrespiSchrage13b}\label{prop:Minty_Min1}
Let $F:X\to \P(Y)$ be a $C$-convex function and $x_0\in \dom F$ a  scalarized weak minimizer, then $x_0$ satisfies
\begin{align}\label{eq:mvi2}
\begin{array}{cl}
&F(x_0)+C=Y\;\vee\\
&\forall x\in X,\, \exists w^*\in W^*:\; \vp^\triup_{F,w^*}(x)\neq-\infty\wedge\of{\vp^\triup_{F,w^*}}^\downarrow(x,x_0-x)\leq 0.
\end{array}
\end{align}
\end{proposition}

Generalizing the above results from $C$-convexity of $F$ to radial $C$-convexity w.r.t $x_0$ is immediate and does not need to be proven.

\begin{corollary}\label{cor:conclusion}
Let $F:X\to\P(Y)$ be radially $C$-convex with $x_0\in\dom F$ and $F(x_0)+C\neq Y$.
Then
\[
\eqref{eq:svi_scalar}\quad\Rightarrow\quad \eqref{eq:w-sc-Min}\quad \Rightarrow\quad \eqref{eq:mvi_scalar}.
\]
If additionally $F(x_0)$ is a compact set, then \eqref{eq:w-sc-Min} is equivalent to \eqref{eq:w-Min} and \eqref{eq:w-l-Min}
and if additionally $F$ is radially upper Hausdorff continuous w.r.t $x_0$  and the scalarizations $\vp^\triup_{F,w^*}$ are proper for all $w^*\in W^*$, then
\[
\eqref{eq:svi_scalar}\quad\Leftrightarrow\quad \eqref{eq:w-Min}\quad \Leftrightarrow\quad \eqref{eq:mvi_scalar}.
\]
\end{corollary}


\section{Appendix}\label{sec:Appendix}
For the readers convenience, we include the proofs of those results quoted from \cite{CrespiHamelSchrage13W},  \cite{CrespiSchrage13b}.

\begin{proposition}\cite[Proposition 2.11]{CrespiSchrage13b}
  Let $A,B\in\P(Y)$ be given such that $B+C$ is convex. Then
  $A\ll B$ implies
  \begin{align*}
  \forall w^*\in W^*:\quad \inf\limits_{a\in A}w^*(a)=-\infty\;\vee\;\inf\limits_{a\in A}w^*(a)<\inf\limits_{b\in B}w^*(b),
  \end{align*}
  which in turn implies $A<B$ if $A+C$ is convex.
\end{proposition}
\proof
Indeed, $A\ll B$ by definition implies $B+U\subseteq A+C$ for some $U\in\mathcal U_Y(0)$, thus
  \begin{align*}
  \forall w^*\in W^*:\inf\limits_{a\in A}w^*(a)\leq \inf\limits_{\substack{b\in B,\\ u\in U}}w^*(b+u)=\inf\limits_{b\in B}w^*(b)+\inf\limits_{u\in U}w^*(u).
  \end{align*}
But as $\inf\limits_{u\in U}w^*(u)<0$ is true for all $w^*\in W^*$, this is the first implication.
As for the second implication, assume $B\nsubseteq \Int (A+C)$ and $A+C$ convex. Then by a separation theorem
  \begin{align*}
  \exists w^*\in W^*:\inf\limits_{b\in B}w^*(b)\leq \inf\limits_{a\in A}w^*(a)\neq-\infty,
  \end{align*}
  as $\Int(A+C)\neq\emptyset$ is assumed.
\pend

\begin{proposition}\cite[Proposition 4.13]{CrespiHamelSchrage13W}
Let $\vp :\R \to \OLR$ be l.s.c. and semistrictly quasiconvex with $\dom \vp\subseteq \sqb{0,1}$. Then there exist $s_0 \leq t_0 \in \sqb{0,1}$ such that $\vp$ is strictly decreasing on $\of{0, s_0}$, strictly increasing on $\of{t_0,1}$ and constantly equal to $\inf\cb{\vp(x)\st x\in X}$ on $\sqb{s_0, t_0}$.
\end{proposition}
\proof
Let $\vp(0)=\vp(1)$ be given, $t\in\sqb{0,1}$. If $\vp(t)>\vp(0)$, then $\vp(s)<\vp(t)$ and thus $\vp(s)=\vp(0)$ is true for all $s\in\sqb{0,1}\setminus\cb{t}$  by semistrict quasiconvexity of $\vp$.
Lower semicontinuity of $\vp$ thus implies $\vp(t)\leq \vp(0)$, a contradiction.
Define the level sets of $\vp$ w.r.t. $t\in\OLR$ as
\begin{align*}
 L^\leq_\vp(t)=\cb{x\in\sqb{0,1}\st \vp(x)\leq t}.
\end{align*}
By the above  $L^\leq_\vp(t)$ is convex and $L^\leq_\vp(t)$ is closed by lower semicontinuity of $\vp$ for all $t\in\OLR$.
Especially,
\[
L^\leq_\vp(\inf\limits_{x\in \sqb{0,1}}\vp(x))=\bigcap\limits_{x\in \sqb{0,1}}L^\leq_\vp(\vp(x))
\]
is a closed convex set, hence either $-\infty$ is attained in some $x\in\sqb{0,1}$, trivially implying $L^\leq_\vp(\inf\limits_{x\in \sqb{0,1}}\vp(x))\neq\emptyset$, or the Weierstrass Theorem implies that the infimum of the loser semicontinuous function $\vp$ is attained on the compact set $\sqb{0,1}$.

Now if $0<s<t\leq s_0$, then semistrict quasiconvexity of $\vp$ implies $\vp(0)>\vp(s)>\vp(s_0)$ and $\vp(s)>\vp(t)>\vp(s_0)$, as $\vp(s_0)=\inf\limits_{x\in \sqb{0,1}}\vp(x)$.
But thus $\vp$ is strictly decreasing on $\sqb{0,s_0}$ and the same arguments prove strict monotonicity on the interval $\sqb{t_0,1}$.
\pend

The following result is Diewert's Mean Value Theorem \cite{Diewert81}.

\begin{proposition}
\label{PropDiewert}
Let $\vp \colon X \to \OLR$ and $a, b\in X$ be such that $\vp_{a,b} \colon \sqb{0,1} \to \R$ is lower semicontinuous and real-valued. Then, there exist $0 \leq t < 1$ and $0 < s \leq 1$ such that
\begin{align*}
 \vp(b) -  \vp(a) &\leq ( \vp_{a,b})^\downarrow(t,1) \; \text{and} \\
 \vp(a) -  \vp(b) &\leq ( \vp_{a,b})^\downarrow(s,-1).
\end{align*}
\end{proposition}

By a careful case study, we can extend this classical result to the case when $\vp_{a,b} \colon \sqb{0,1}\to \OLR$ is extended real-valued and not necessarily proper. Then, the difference has to be replaced by the inf-residual in $\OLR$,
\[
\forall s,t\in\OLR:\quad s\idif t=\inf\cb{r\in\R\st s\leq t+r},
\]
assuming $(+\infty)+r=+\infty$ and $(-\infty)+r=-\infty$ for all $r\in \R$.
Especially
\[
\forall s\in\OLR:\quad (-\infty)\idif s=s\idif (+\infty)=-\infty.
\]

\begin{theorem}\label{thm:Diewert}\cite[Proposition 4.2]{CrespiHamelSchrage13W}
Let $\vp \colon X \to \OLR$ and $a, b \in X$ be given such that $a \neq b$ and $\vp_{a,b} \colon \R \to \OLR$ is lower semicontinuous. Then:

(a) If either $\vp(a) = +\infty$, or $\cb{a,b} \subseteq \dom \vp$, then there exists $0\leq t< 1$ such that
\[
  \vp(b)\idif \vp(a)\leq \of{\vp_{a,b}}^\downarrow(t,1).
\]

(b) If either $\vp(b)=+\infty$, or $\cb{a,b} \subseteq \dom \vp$, then there exists $0< s \leq 1$ such that
\[
\vp(a) \idif \vp(b) \leq \of{\vp_{a,b}}^\downarrow(s,-1).
\]
\end{theorem}

{\sc Proof.} (a) The proof of the first inequality is given via a case study. If $\vp(a)=+\infty$ or $\vp(b)=-\infty$, then
\[
\vp(b)\idif\vp(a)=\inf\cb{r\in \R\st \vp(b)\leq \vp(a)+r}=-\infty,
\]
so the first inequality is trivially satisfied.

Next, assume $\cb{a,b}\subseteq\dom\vp$ and $\vp(b)\neq-\infty$. If $\vp_{a,b}(t)=-\infty$ for some $0\leq t<1$, then by lower semicontinuity $\vp_{a,b}(t_0)=-\infty$, setting
\[
  t_0=\sup\cb{t\in\cb{0,1}\st \vp_{a,b}(t)=-\infty}
\]
and by assumption $t_0<1$. Hence $\of{\vp_{a,b}}^\downarrow(t_0,1)=+\infty$, satisfying the first inequality.

Finally, let $\cb{a,b}\subseteq\dom\vp$ and $\vp(b)\neq-\infty$ be assumed and $\vp_{a,b}(t)=+\infty$ for some $0<t<1$ and set
\[
  t_0=\inf\cb{t\in\of{0,1}\st \vp_{a,b}(t)=+\infty}.
\]
If $t_0=0$, then we are finished, as in this case $\of{\vp_{a,b}}^\downarrow(0,1)=+\infty$ is true, hence assume $0<t_0$. In this case, $\sqb{0,t}\subseteq\dom\vp_{a,b}$ is true for all $t\in\of{0,t_0}$, and the above result combined with Proposition \ref{PropDiewert} applied to $b=a+t(b-a)$ gives that for all $0 < t < t_0$ there exists a $0\leq \bar t<1$ such that
\[
  \vp(a+t(b-a))\leq \vp(a)\isum \of{\vp_{a,a+t(b-a)}}^\downarrow(\bar t,1),
\]
But as $\of{\vp_{a,a+t(b-a)}}^\downarrow(\bar t,1)=\of{\vp_{a,b}}^\downarrow(\bar t,1)$ is true and by lower semicontinuity of $\vp_{a,b}$ the value $\vp(a+t(b-a))$ converges to $+\infty$ as $t$ converges to $t_0$, this implies that $\of{\vp_{a,b}}^\downarrow(\bar t,1)$ converges to $+\infty$ and eventually satisfies the desired inequality.

(b) Notice that $\vp_{a,b}(s)=\vp_{b,a}(1-s)$ and $\of{\vp_{a,b}}^\downarrow(s,-1)=\of{\vp_{b,a}}^\downarrow((1-s),1)$, hence the result is immediate from the above. \pend

\begin{proposition}\cite[Proposition 4.14]{CrespiHamelSchrage13W}
If $\dom \vp$ is star shaped at $x_0$ and $\vp$ is radially pseudoconvex and l.s.c. w.r.t. $x_0$, then it is radially semistrictly quasiconvex w.r.t. $x_0$.
\end{proposition}
\proof Assume that for some $b\in\dom \vp$ the function $\vp_{a,b}$ is not semistrictly quasiconvex. Then there are $r, s, t \in \R$ such that
$0 \leq r < s < t \leq 1$, $\vp_{a,b}\of{r} \neq \vp_{a,b}(t)$ and
\[
\max\cb{\vp_{a,b}\of{r}, \vp_{a,b}(t)} \leq \vp_{a,b}(s).
\]
We assume $\vp_{a,b}\of{r} < \max\cb{\vp_{a,b}\of{r}, \vp_{a,b}\of{t}} = \vp_{a,b}\of{t}$. The other case can be dealt with by symmetric arguments.

Fix $\delta > 0$ such that $\vp_{a,b}\of{r} < \vp_{a,b}\of{t} - \delta$. Since $\vp_{a,b}$ is l.s.c. the set
\[
\cb{s' \in \R \mid \vp_{a,b}\of{s'} > \vp_{a,b}\of{t} - \delta}
\]
is open. Hence there is $\eps > 0$ such that $[s - \eps, s + \eps] \subseteq \of{r, t}$ and
\[
\forall s' \in [s - \eps, s + \eps] \colon   \vp_{a,b}\of{t}  - \delta < \vp_{a,b}\of{s'} \in\R.
\]

Take $s' \in [s, s+\eps)$, $s'' \in (s', s+\eps]$ and assume $\vp_{a,b}\of{s''} < \vp_{a,b}\of{s'}$. By Diewerts Mean-Value-Theorem \ref{PropDiewert} there exists an $\hat s \in (s' , s'']$ satisfying
\[
0 < \vp_{a,b}\of{s'} - \vp_{a,b}\of{s''} \leq \of{\vp_{a,b}}^\downarrow\of{\hat s, s' - s''}.
\]
Indeed, setting $a' = a+s'(b-a)$, $b' = a+s''(b-a)$ one obtains by Diewerts Mean-Value-Theorem an $\alpha \in (0, 1]$ satisfying $\vp\of{a'} - \vp\of{b'} \leq \of{\vp_{a,b}}^\downarrow\of{\alpha, -1}$. Defining $\hat s = s + \alpha(s''-s') \in (s', s'']$ and observing $\vp\of{a'} = \vp_{a,b}\of{s'}$, $\vp\of{b'} = \vp_{a,b}\of{s''}$ and $\of{\vp_{a,b}}^\downarrow\of{\alpha, -1} = \of{\vp_{a,b}}^\downarrow\of{\hat s, s' - s''}$ one obtains the above inequality. Using the positive homogeneity of the directional derivative we can multiply the inequality $0 < \of{\vp_{a,b}}^\downarrow\of{\hat s, s' - s''}$ by $\frac{r - \hat s}{s' - s''} > 0$ and obtain
$0 < \of{\vp_{a,b}}^\downarrow\of{\hat s, r - \hat s}$. The pseudoconvexity of $\vp_{a,b}$ yields $\vp_{a,b}\of{r} \geq \vp_{a,b}\of{\hat s}$ which contradicts the assumption $\vp_{a,b}\of{r} < \vp_{a,b}\of{t} - \delta < \vp_{a,b}\of{\hat s} - \delta$ (observe $\hat s \in [s, s+\eps]$). Hence $\vp_{a,b}\of{s''} \geq \vp_{a,b}\of{s'}$ whenever $s', s'' \in [s, s+\eps]$ and $s' < s''$. This implies
\[
\forall s' \in [s, s+\eps) \colon \of{\vp_{a,b}}^\downarrow\of{s', 1} \geq 0,
\]
and positive homogeneity of the directional derivative implies
\[
\of{\vp_{a,b}}^\downarrow\of{s', t-s'} \geq 0
\]
 and this by pseudoconvexity of $\vp_{a,b}$
\[
\vp_{a,b}\of{t} \geq \vp_{a,b}\of{s'} \geq \vp_{a,b}\of{s} \geq \vp_{a,b}\of{t}.
\]
This means $\vp_{a,b}\of{s'} = \vp_{a,b}\of{t}$ for all $s' \in [s, s+\eps)$. In turn, this implies that for $s' \in (s, s+\eps)$ we have $\of{\vp_{a,b}}^\downarrow\of{s', -1} \geq 0$, hence $\of{\vp_{a,b}}^\downarrow\of{s', r- s'} \geq 0$ and by pseudoconvexity $\vp_{a,b}\of{s'} \leq \vp_{a,b}\of{r}$. This contradicts the assumption $\vp_{a,b}\of{r} < \vp_{a,b}\of{t}$, hence (together with the symmetric case) the function $\vp_{a,b}$ is semistrictly quasiconvex for all $b \in \dom \vp$.\pend

\begin{proposition}\cite[Proposition 4.5]{CrespiSchrage13b}
Let $F:X\to \P(Y)$ be a $C$-convex function, $x_0\in \dom F$.
If $x_0$ solves the scalarized Stampacchia variational inequality \eqref{eq:svi2App}, then it is a scalarized weak minimizer.
\begin{align}\label{eq:svi2App}
F(x_0)+C=Y\quad \vee\quad\forall x\in X\, \exists w^*\in W^*:\quad \of{\vp^\triup_{F,w^*}}^\downarrow(x_0,x-x_0)\geq 0
\end{align}
\end{proposition}
\proof
Assume to the contrary that $F(x_0)+C\neq Y$ and it exists $x\in X$ such that
\begin{align*}
\forall w^*\in W^*:\quad \vp^\triup_{F,w^*}(x)<\vp^\triup_{F,w^*}(x_0).
\end{align*}
As all scalarizations are convex by assumption, this contradicts \eqref{eq:svi2App}.
\pend

\begin{proposition}\cite[Lemma 4.9]{CrespiSchrage13b}
Let $F:X\to \P(Y)$ be a $C$-convex function, $x_0\in \dom F$.
 If $x_0$ satisfies \eqref{eq:mvi2_App} , then it is a scalarized weak minimizer.
\begin{align}\label{eq:mvi2_App}
\begin{array}{cl}
&F(x_0)+C=Y\;\vee\\
&\forall x\in X\, \exists w^*\in W^*:\quad \vp^\triup_{F,w^*}(x)\neq-\infty\wedge\of{\vp^\triup_{F,w^*}}^\downarrow(x,x_0-x)\leq 0.
\end{array}
\end{align}
\end{proposition}
\proof
If $x_0$ is a scalarized weak minimizer, then either $F(x_0)+C=Y$ or for every $x\in X$ there exists a $w^*\in W^*$ such that $\vp^\triup_{F,w^*}(x_0)\leq \vp^\triup_{F,w^*}(x)\neq-\infty$ and thus
\[
\of{\vp^\triup_{F,w^*}}^\downarrow(x,x_0-x)\leq
\vp^\triup_{F,w^*}(x_0)- \vp^\triup_{F,w^*}(x)\leq 0.
\]
\pend


\end{document}